\documentclass[12pt]{amsart}
\usepackage{amscd,amsthm,amsfonts,amssymb,amsmath,euscript}
\usepackage[dvips]{graphicx}
\usepackage[dvips]{graphics}
\usepackage[T2A]{fontenc}
\usepackage[cp866]{inputenc}

\newtheorem{theorem}{Theorem}

\newcommand{\Aut}{\mathop{\sf Aut}\nolimits}

\newcommand{\ord}{\mathop{\sf ord}\nolimits}

\begin{document}

\centerline{\bf Test for reality of algebraic functions}

{\ }

\centerline{\bf S.M.Natanzon{\footnote{Supported by grants
RFBR-04-01-00762, NWO 047.011.2004.026 (RFBR 05-02-89000-НВО-а),
NSh-4719.2006.1, MTM 2005-01637, CRDF RM1-2543-MO-03, INTAS
05-7805}} }

\vskip 0.6cm

In this paper I prove the next affirmation
\begin{theorem}\label{t1} 1) A complex algebraic function on a real
algebraic curve with real points is equivalent to a real algebraic
function, if and only if the divisor of preimage of its critical
values is stable under the involution of complex conjugation. 2) A
complex algebraic function on a real algebraic curve without real
points is equivalent to a real or pseudoreal algebraic function, if
and only if the divisor of preimage of its critical values is stable
under the involution of complex conjugation.
\end{theorem}

The proof uses representations of  algebraic functions by Fuchsian
groups as [N1, N(sect.1.6)]. It is interesting to compare this
theorem with a conjecture of B. and M. Shapiro. It states, that a
complex algebraic function is equivalent to real, if all its
critical points are real. Proofs of special cases of this
conjecture, when genus and degree satisfy $\deg(f)^2-4\deg(f)+3<3g$
or $\deg(f)\leq 4$, are results of papers [EG1,EG2, ESS, DEISS].

{\ }

Let us express our theorem in term of Riemann surfaces. A real
algebraic curve is equivalent to a pair $(P,\tau)$, where $P$ is a
compact Riemann surface (the set of complex points of the curve) and
$\tau: P\rightarrow P$ is an antiholomorphic involution (the
involution of complex conjugation). Fixed points of $\tau$ are real
points of the real algebraic curve[AG].

A complex algebraic function on a real curve $(P,\tau)$ is a
holomorphic map
$f:P\rightarrow\hat{\mathbb{C}}=\mathbb{C}\cup\infty$. The number
$\ord_p$ is called the number of ramification of $f$ in $p\in P$, if
$f(z)=z^{\ord_p}$ in some local map  $z:U\rightarrow\mathbb{C}$,
where $z(p)=0$. The points $\Sigma_P=\{p\in P|\ord_p>1\}$ form the
set of critical points. Its image
$\Sigma_{\hat{\mathbb{C}}}=f(\Sigma_P)\subset\hat{\mathbb{C}}$ form
the set of critical values. Let us consider the preimage of critical
values $\Sigma_{P\hat{\mathbb{C}}}=
f^{-1}(\Sigma_{\hat{\mathbb{C}}})\subset P$ and its divisor
$\Sigma(f)=\sum_{p\in \Sigma}\ord_pp$

Functions $f$ and $\tilde{f}$ are called equivalent if
$\tilde{f}=gf$ for $g\in \Aut(\hat{\mathbb{C}})$. A complex
algebraic function is called real, if $\overline{f(\tau(p))}=f(p)$,
and pseudoreal if $\overline{f(\tau(p))}=-\frac{1}{f(p)}$ [N2,
N(sect.3.4)].

Let us illustrate this definitions by example. Consider a polynomial
$F(x,y)$ with real coefficients. Then the affine real algebraic
curve $F(x,y)=0$ generates the pair $(P_F,\tau_F)$, where $P_F$ is a
regularization for complex points
$\{(x,y)\in\mathbb{C}^2|F(x,y)=0\}$ of the curve and
$\tau_F(x,y)=(\bar{x},\bar{y})$. The set of real points of the real
curve is the closure of $\{(x,y)\in\mathbb{R}^2|F(x,y)=0\}$. Complex
algebraic functions are generated by restrictions of complex
polynomials $G(x,y)$ on $\{(x,y)\in\mathbb{C}^2|F(x,y)=0\}$. The
function is real if and only if the polynomial $G(x,y)$ has real
coefficients.

Thus, theorem \ref{t1} is equivalent to

\begin{theorem} Let $P$ be a compact Riemann surface,
$f:P\rightarrow\hat{\mathbb{C}}$ be a holomorphic map, $\tau:
P\rightarrow P$ be a antiholomorphic involution. Then $\tau\Sigma(f)
=\Sigma(f)$, if an only if $\overline{gf\tau}=gf$ or
$\overline{gf\tau}=-\frac{1}{gf}$ for $g\in \Aut(\hat{\mathbb{C}})$,
and, moreover, $\overline{gf\tau}=gf$ for $\{p\in P|\tau
p=p\}\neq\emptyset$.
\end{theorem}

\textsl{Proof}. It is obviously that $\tau\Sigma(f) =\Sigma(f)$
follows from $\overline{gf\tau}=gf$ and
$\overline{gf\tau}=-\frac{1}{gf}$. We use the induction by the
degree $\deg(f)$ in order that to prove that $\tau\Sigma(f)
=\Sigma(f)$ generate $\overline{gf\tau}=gf$ or
$\overline{gf\tau}=-\frac{1}{gf}$. The theorem is obviously for
$\deg(f)=0$. Let $\deg(f)>0$ and, thus, $\Sigma(f)\neq0$. Consider
an uniformisation $\Lambda\rightarrow S$ of
$S=\hat{\mathbb{C}}\setminus f(\Sigma(f))$ by a Fuchsian group
$\Gamma$ on $\Lambda=\{z\in\mathbb{C}|\sf Im(z)>0\}$. Then a
subgroup $\tilde{\Gamma}\subset\Gamma$ uniformises the Riemann
surface $\tilde{P}=P\setminus f^{-1}(f(\Sigma(f)))$. The involution
$\tau$ is generated by an antiholomorphic map
$\check{\tau}:\Lambda\rightarrow\Lambda$, where
$\check{\tau}\tilde{\Gamma}\check{\tau}^{-1}=\tilde{\Gamma}$.

Consider $x\in \Sigma(f)$. Let $\check{x}\in\mathbb{R}$ be the fixed
point of a parabolic automorphism $\tilde{\gamma}\in\tilde{\Gamma}$,
that correspond to puncture $x$. Then the stationary subgroup
$\tilde{\Gamma}_{\check{x}}=\{\gamma\in
\tilde{\Gamma}|\gamma(\check{x})=\check{x}\}$ is a subgroup of index
$\ord_x(f)>1$ of $\Gamma_{\check{x}}=\{\gamma\in
\Gamma|\gamma(\check{x})=\check{x}\}$. Thus, using $\tau\Sigma(f)
=\Sigma(f)$ we have
$\check{\tau}\tilde{\Gamma}_{\check{x}}\check{\tau}^{-1}\subset\tilde{\Gamma}$
and $\check{\tau}\Gamma_{\check{x}}\check{\tau}^{-1}\subset\Gamma$.
Moreover, the condition $\tau\Sigma(f) =\Sigma(f)$ is equivalent to
coincidence of indexes of stationer subgroups
$[\Gamma_{\tilde{x}}:\tilde{\Gamma}_{\tilde{x}}]=
[\Gamma_{\check{\tau}(\tilde{x})}:\tilde{\Gamma}_{\check{\tau}(\tilde{x})}]$
for all fixed points $\tilde{x}$ of parabolic automorphisms from
$\tilde{\Gamma}$.

Let $\hat{\Gamma}\subset\Gamma$ be the subgroup, generated by
$\tilde{\Gamma}$, $\Gamma_{\tilde{x}}$ and
$\check{\tau}\Gamma_{\check{x}}\check{\tau}^{-1}$. It contains
$\tilde{\Gamma}$ as a proper subgroup. Moreover,
$\check{\tau}\hat{\Gamma}\check{\tau}^{-1}=\hat{\Gamma}$ and
$[\Gamma_{\tilde{x}}:\hat{\Gamma}_{\tilde{x}}]=
[\Gamma_{\hat{\tau}(\tilde{x})}:\hat{\Gamma}_{\hat{\tau}(\tilde{x})}]$
for all fixed points $\tilde{x}$ of parabolic automorphisms from
$\hat{\Gamma}$.

Consider $\hat{P}=\Lambda/\hat{\Gamma}$. The inclusions of groups
$\tilde{\Gamma}\subset\hat{\Gamma}\subset\Gamma$ generate the
coverings
$\Lambda\to\tilde{P}\overset{\tilde{\psi}}\to\hat{P}\overset{\hat{\psi}}\to
S$.  These coverings generate holomorphic maps of compact surfaces
$P\overset\varphi\to\dot{P}\overset{\dot{\varphi}}\to\hat{\mathbb{C}}$
after gluing of punctures. Moreover, $f=\dot{\varphi}\varphi$, and
$\deg(\dot{\varphi})<\deg(f)$. The involution $\check{\tau}$
generates an antiholomorphic involution $\dot{\tau}$ on $\dot{P}$.
Moreover, $\varphi\tau=\dot{\tau}\varphi$, and
$\dot{\tau}\Sigma(\dot{\varphi}) =\Sigma(\dot{\varphi})$. Therefore,
by inductive hypothesis $\overline{\dot{\varphi}\dot{\tau}}$ is
equal to $g\dot{\varphi}$  or $-\frac{1}{g\dot{\varphi}}$ for $g\in
\Aut(\hat{\mathbb{C}})$. Thus, the function
$\overline{f\tau}=\overline{\dot{\varphi}\varphi\tau}=
\overline{\dot{\varphi}\dot{\tau}\varphi}$ is equal
$g\dot{\varphi}\varphi=gf$ or
$-\frac{1}{g\dot{\varphi}\varphi}=-\frac{1}{gf}$. It is follow from
[N2,N(sect.3.4)], that $\{p\in P|\tau p=p\}=\emptyset$ in the last
case.

$\Box$

\centerline{References}

[AG] Alling N.L., Greenleaf N., Foundations of the theory of Klein
surfaces. Berlin - Heidelbeс rg - N.Y.: Springer - Verlag, 1971.
Lecture Notes in Mathematics. V.219.

[EG1]  Eremenko A., Gabrielov A., Elementary proof of the B. and M.
Shapiro conjecture for rational functions, 2005, math/0512370v1, 21
p.

[EG2]  Eremenko A., Gabrielov A., Rational functions with real
critical points and the B. and M. Shapiro conjecture in real
enumerative geometry. Ann. of Math. (2) 155 (2002), no.1, 105-129.

[ESS] Ekedahl T, Shapiro B., Shapiro M., First steps towards total
reality of meromorphic functions, Mosc. Math. J. 6 (2006), no.1,
95-106.

[DEISS] Degtyarev A., Ekedahl T., Itenberg I., Shapiro B., Shapiro
M., On total reality of meromorphic functions, Ann. Inst. Fourier
(Grenoble) 57 (2007), no. 6, 2015-2030 ( math/0605077).

[N1] Natanzon S.M., Uniformization of spaces of meromorphic
functions. Dokl. Acad. Nauk SSSR, 287:5 (1986), 1058-1061 (Russian);
English: Soviet. Math. Dokl., 33:2(1986), 487-490.

[N2] Natanzon S.M., Topology of 2-dimensional coverings and
meromorphic functions on real and complex algebraic curves. Selecta
Math.Soviet., 12:3 (1993), 251-291.

[N]  Natanzon S.M., Moduli of Riemann surfaces, real algebraic
curves, and their superanalogs. Translations of Mathematical
Monographs, American Mathematical Society, Vol.225, (2004), 160 P.

\end{document}